\documentclass[12pt]{amsart}
\newtheorem{theorem}{Theorem}[section]
\newtheorem{lemma}[theorem]{Lemma}
\newtheorem{corollary}[theorem]{Corollary}
\theoremstyle{definition}

\theoremstyle{remark}
\newtheorem{remark}[theorem]{Remark}
\numberwithin{equation}{section}

\newsymbol\subsetneq 2328
\newsymbol\dotplus 1275

\def\z*{\bar z}
\def\B{\mathsf B}
\def\GM{\mathcal G}
\def\H{\mathcal H}
\def\fh{\mathfrak h}
\def\HP{\H_{+}}
\def\K{\mathcal K}
\def\N{\mathcal N}
\def\RE{\mathbb R}
\def\C{{\mathbb C}}
\def\trip{\left\{\fh,\gamma_1,\gamma_2\right\}}

\def\G*{G_*}
\def\ph*{\phi_*}

\begin{document}

\title
[Boundary Triples for Singular Perturbations]
{Boundary Triples and Weyl Functions for Singular Perturbations
of Self-Adjoint Operators\footnote{To appear in: {\it Methods of Functional Analysis and Topology}}}
\author{Andrea Posilicano}

\address{Dipartimento di Scienze, Universit\`a dell'Insubria, I-22100
Como, Italy}

\email{posilicano@uninsubria.it}

\begin{abstract}
Given the symmetric operator $A_\N$ obtained by restricting the
self-adjoint operator $A$ to $\N$, a linear dense set, closed with 
respect to the graph norm, we determine a convenient boundary triple for the
adjoint $A_N^*$ and the corresponding Weyl function. These objects provide us 
with the self-adjoint extensions of $A_\N$ and their resolvents.
\end{abstract}

\maketitle

\section{Introduction}
Let $A:D(A)\subseteq\H\to\H$ be a self-adjoint operator on the
Hilbert space $\H$. Another self-adjoint operator $\widehat A$ is said to
be a singular perturbation of $A$ if the set $\N:=\{\phi\in
D(\widehat A)\cap D(A)\, :\, \widehat A\phi=A\phi\}$ is dense in
$\H$ (see e.g. \cite{[K]}). Since $\N$ is closed with respect to the graph norm on $D(A)$, 
the linear operator $A_\N$, obtained by restricting $A$ to $\N$, is a
densely defined closed symmetric operator. Thus $\widehat A$ is a
singular perturbation of $A$ if and only if it is a self-adjoint
extension of $A_\N$ such that $D(\widehat A)\cap D(A)=\N$, where
$\N\subsetneq D(A)$
is any dense set which is closed with respect to the graph
norm on $D(A)$. Therefore all singular perturbations of $A$ could be 
determined by using von
Neumann's theory \cite{[N]}. By such a theory, given a closed densely
defined symmetric operator $S$, one has
$$
D(S^*)= D(S)\oplus
\K_+\oplus \K_-\,,\quad
S^*(\phi_0+\phi_++\phi_-)=S\phi_0+i\phi_+-i\phi_-\,,
$$ 
where the direct sum decomposition is orthogonal with
respect to the graph inner product of $S^*$ and 
$\K_\pm:=$Ker$(-S^*\pm i)$ denotes the deficiency
spaces. Any self-adjoint
extension $A_U$ of $S$ is then obtained
by restricting $S^*$ to a subspace of the kind $D(S)\oplus \text{\rm
Graph}\,U$, where
$U:\K_+\to \K_-$ is unitary. Alternatively one could determine the
singular perturbations of $A$ by Kre\u\i n's 
resolvent formula (see \cite{[K1]}, \cite{[K2]},  
\cite{[S]} for the cases where dim$\,\K_\pm=1$, dim$\,\K_\pm<+\infty$,
dim$\,\K_\pm=+\infty$ respectively; also see \cite{[GMT]}). This
approach gives the resolvent difference of any pair (in our situation 
$A$ and $\hat A$) of self-adjoint 
extensions of $S$, 
thus allowing for a better understanding of the spectral
properties of such extensions. 
Like in von Neumann's theory, also in Kre\u\i n's one a main
role is played by the defect spaces, through the orthogonal projections
onto $\K_\pm$. This is a conseguence of the fact that both 
these theories regard self-adjoint
extensions of an arbitrary closed densely defined symmetric operator
with equal defect indices. However in the case of singular perturbation 
the situation is simpler. Indeed here $A_\N$ is not an arbitrary
symmetric operator but is the restriction to $\N$ of a given
self-adjoint operator $A$. Thus, instead
of the orthogonal projections
onto $\K_\pm$, one can use the orthogonal projection $\pi\, :\,
\HP\to\N^\perp$, where $\HP$
denotes the Hilbert space given by the set
$D(A)$ equipped with the scalar product $\langle\phi_1,\phi_2\rangle_{+}
:=\langle A\phi_1,A\phi_2\rangle+\langle\phi_1,\phi_2\rangle$, 
and the orthogonal decomposition $\HP=\N\oplus\N^\perp$ is used,
being $\N$ closed in $\HP$. More generally, 
since this gives advantages in concrete applications where
usually a variant of $\pi$ is what is known in advance, one can
consider a bounded linear map $\tau:\HP\to\fh$, $\fh$ an auxiliary
Hilbert space, such that $\N:=$Ker$(\tau)$ is
dense in $\H$ and Ran$(\tau)=\fh$, so that  
$\HP\simeq \text{\rm Ker$(\tau)$}\oplus\fh$. This alternative approach has been developed in
\cite{[P1]} as regards Kre\u\i n's formula and in \cite{[P3]}, where an
additive decomposition of any singular perturbation is given and the
explicit connection with von Neumann's theory is found. 
The approach contained in
\cite{[P1]}, \cite{[P2]}, \cite{[P3]}, 
looks simpler than the original ones (for example no
knowledge of either $A_\N^*$ or $\K_\pm$ is needed), allows for a
natural formulation in terms of (abstract)
boundary conditions and makes easier to work
out concrete applications where $\tau$ is the trace (restriction) map
along some null subset of $\RE^d$ and $A$ is a (pseudo-)differential
operator (see the examples contained in the quoted references). \par
Another approach, different from both von Neumann's and Kre\u\i n's
theories, has been used to obtain
self-adjoint extensions of a given symmetric operator $S$ with equal
defect indices: this is the theory of boundary triplets introduced by 
Bruk and Kochubei \cite{[B]}, \cite{[Koch]} and then 
successively developed in many papers and books 
(see e.g. \cite{[GG]}, \cite{[DM1]}, \cite{[DM2]} and references therein). 
Here one needs to find a boundary triple $\trip$ for $S^*$, i.e. one
needs a Hilbert space $\fh$ (with scalar product 
$[\cdot,\cdot]$) and two linear maps 
$\gamma_1$ and $\gamma_2$ on $D(S^*)$ to $\fh$
such that $\phi\mapsto (\gamma_1\phi,\gamma_2\phi)$ is surjective and
$$
\langle S^*\phi,\psi\rangle-\langle \phi,S^*\psi\rangle=
[\gamma_1\,\phi,\gamma_2\,\psi]-
[\gamma_2\,\phi,\gamma_1\,\psi]\,.
$$
Once a boundary triple is known, 
any self-adjoint extensions of $S$ is then obtained by
restricting $S^*$ to the set of $\phi$'s such that the couple
$(\gamma_2\phi,\gamma_1\phi)$ belongs to a self-adjoint relation. \par
Boundary triples theory, thanks to the
concept of Weyl function successively introduced in \cite{[DM3]}, generalizing 
a concept earlier used by Weyl in the study of Sturm-Liouville
problems, allows for a spectral analysis of the self-adjoint
extensions of $S$ (see e.g. \cite{[DM1]}, \cite{[BMN]} and references 
therein).  \par
The scope of this paper is to work out the theory of boundary triples and
Weyl functions in
the case of singular perturbations, making use of the orthogonal projection
$\pi\,:\,\HP\to\N^\perp$ or better of its generalization given by a map 
$\tau:\HP\to\fh$. This is done in section 3, after that, 
in section 2, we have given a concise review (we refers to 
\cite{[DM1]}, \cite{[GG]} for the proofs) of
the theory of boundary triples. In particular a convenient boundary
triple for $A^*_\N$, and
the corresponding Weyl function are given in Theorem 3.1. The
successive Corollary 3.2 characterizes all singular perturbations of
$A$ in terms of (abstract) boundary conditions of the kind
$\Theta\zeta_\phi=\tau\phi_*$, where $\Theta$ is self-adjoint on $\fh$. 
As we already said, in the case $A$ is a differential 
operator, these are indeed concrete boundary
conditions since usually $\tau$ is the (trace) evaluation map along some null
subset. Finally in Theorem
3.4 we determine the possible eigenvectors (and their multiplicity) of
the singular perturbations of a self-adjoint operator. 

\section{Boundary triples and Weyl Functions.}
Let $S: D(S)\subseteq\H\to\H$ be a densely defined closed symmetric
operator on the Hilbert space $\H$ with inner product 
$\langle\cdot,\cdot\rangle$ and corresponding norm $\|\cdot\|$. We
will suppose that $S$ has equal defect indices $n_+=n_-$, where 
$n_\pm:=$dim$\,\K_\pm$ and $\K_\pm:=$Ker$(-S^*\pm i)$. \par 
A triple 
$\trip$, where $\fh$ is a Hilbert space with inner product 
$[\cdot,\cdot]$ and 
$$\gamma_1\, :\, D(S^*)\to\fh\,,\quad\gamma_2\, :\, D(S^*)\to\fh\,,$$ 
are two linear maps such that 
$$\gamma:D(S^*)\to\fh\oplus\fh\,,\quad \gamma\phi:=(\gamma_1\phi,\gamma_2\phi)$$
is surjective
and
$$
\langle S^*\phi,\psi\rangle-\langle \phi,S^*\psi\rangle=
[\gamma_1\,\phi,\gamma_2\,\psi]-
[\gamma_2\,\phi,\gamma_1\,\psi]\,.
$$
The definition of boundary triple is well posed. Indeed let 
$$P_\pm\,:\,D(S^*)\to\K_\pm$$
denotes the orthogonal projection given the decomposition 
$$D(S^*)= D(S)\oplus
\K_+\oplus \K_-\,,$$ where the direct sum is orthogonal with
respect to the graph inner product of $S^*$. Then 
$$\gamma_1:=i\,P_+-i\,U P_-\,,\qquad\gamma_2:=P_++U P_-\,,$$  
where $\fh=\K_+$ and $U:\K_-\to\K_+$ is an isometry, give a boundary
triple for $S^*$. This also
shows, since $U$ is arbitrary, that a boundary triple 
is not unique. 
\par 
A closed subspace $\GM\subset \fh\oplus\fh$ is said to be a symmetric
closed relation if 
$$
\forall\,(f_1,g_1),\,(f_2,g_2)\in\GM\,\qquad[f_1,g_2]=
[g_1,f_2 ]\,.
$$
$\GM$ is then said to be a self-adjoint relation if it is maximal
symmetric, i.e. if it does not exists a closed symmetric relation
$\widehat \GM$ such that $\GM\subsetneq\widehat\GM$. Of course the
graph of a self-adjoint operator is a particular case of self-adjoint
relation. \par 
The main result of boundary triples theory is given in the
following
\begin{theorem} \text{\rm (\cite{[GG]}, Theorem 1.6, chapter 3)} Let $\trip$ be a boundary triple
for $S^*$. Then any self-adjoint extension of $S$ is of the kind
$S^*_\GM$, where $S^*_\GM$ denotes
the restriction of $S^*$ to the subspace 
$\{\phi\in D(S^*) : (\gamma_2\phi,\gamma_1\phi)\in\GM\}$, 
$\GM$ being a self-adjoint relation.
\end{theorem} 
Form now on we will denote by $A_\Theta$ the self-adjoint extension which
corresponds to $\GM=$graph$(-\Theta)$, where $\Theta$ is a
self-adjoint operator on $\fh$, and by $A$ the self-adjoint extension 
corresponding to the self-adjoint relation $\GM=\left\{0\right\}\times\fh$.\par Given the boundary triple $\trip$ for $S^*$, the Weyl
function of $S$ corresponding to $\trip$ is defined as the unique map 
$$\Gamma\,:\,\rho(A)\to\B(\fh)$$ such that 
$$
\forall\,\phi_z\in\K_z:=\text{\rm Ker}(-S^*+z)\,,\qquad
\Gamma(z)\gamma_2\, \phi_z=\gamma_1\, \phi_z\,.
$$
Since, for any $z\in \rho(A)$, $\gamma_1$ and $\gamma_2$ are
bijections on $\K_z$ to $\fh$, one can define
$$
G(z):=\left({\gamma_2}_{|\K_z}\right)^{-1}\,,
$$
and thus
$$
\Gamma(z)=\gamma_1 G(z)\,.
$$
Moreover
$$
\Gamma(z)-\Gamma(w)^*=(z-\bar w)\,G(w)^* G(z)\,,
$$
i.e. $\Gamma$ is a $Q$-function of $S$ belonging to the extension $A$
in the sense of Kre\u\i n. 
\par
Weyl functions can be used to deduce spectral properties of the extensions:
\begin{theorem} \text{\rm (\cite{[DM1]}, 
Propositions 1 and 2, section 2)}
$$
z\in\rho(A_\Theta)\cap\rho(A)\iff 
0\in\rho(\Theta+\Gamma(z))
$$
and
$$
\lambda\in\sigma_i(A_\Theta)\cap\rho(A)\iff  0\in\sigma_i(\Theta+\Gamma(\lambda))\,,\quad i=p,c,r\,,
$$
where $\sigma_p(A_\Theta)$, $\sigma_c(A_\Theta)$, $\sigma_r(A_\Theta)$
denote the point, continuous and residual spectrun
respectively. Moreover one has the Kre\u\i n's formula
$$
(-A_\Theta+z)^{-1}=(-A+z)^{-1}+G(z)(\Theta+\Gamma(z))^{-1}G(\bar z)^*\,.
$$
\end{theorem}
The
Weyl function $\Gamma$ is a Herglotz or Nevanlinna operator-valued
function. This means that $\Gamma$ is holomorphic in the upper half
complex plane $\C_+$
and the operator $\Gamma(z)$ is dissipative, i.e. 
$$
\forall\,z\in\C_+\,,\qquad 
\frac{1}{2i}\,\left(\Gamma(z)-\Gamma(z)^*\right)\ge 0\,.
$$
Thus, by the celebrated Nevanlinna-Riesz-Herglotz decomposition, one has 
$$
\Gamma(z)=C+\int_R d\Sigma(t)\,\left(\frac{1}{t-z}-\frac{t}{1+t^2}\right)\,,
$$
where $C\in\B(\fh)$ is self-adjoint and the $\B(\fh)$-valued measure $\Sigma$
is self-adjont and such that
$$
\int_R\frac{d\Sigma(t)}{1+t^2}\ \in\B(\fh)\,.
$$
Given a (unbounded) self-adjoint operator $\Theta$ on $\fh$ and a
$\B(\fh)$-valued Herglotz function $\Gamma\,:\,\C_+\to\B(\fh)$, 
$\Theta+\Gamma(z)$ is boundedly invertible for any $z\in\C_+$ (see
\cite{[P1]}, Proposition 2.1) and $\left(\Theta+\Gamma\right)^{-1}$ is
again a $\B(\fh)$-valued Herglotz function. 
\begin{lemma} \text{\rm (\cite{[BMN]}, Lemma 3.2)} 
The spectral measure of $A_\Theta$ is equivalent 
to the self-adjoint measure appearing in the Nevanlinna-Riesz-Herglotz 
decomposition of $\left(\Theta+\Gamma\right)^{-1}$. 
\end{lemma}  
\section{singular perturbations}
Let $A:D(A)\subseteq\H\to\H$ be a
self-adjoint operator on the Hilbert space
$(\,\H,\langle\cdot,\cdot\rangle\,)$. 
We denote by $(\,\HP, \langle\cdot,\cdot\rangle_+)$ 
the Hilbert space given by the set
$D(A)$ equipped with the scalar product $\langle\cdot,\cdot\rangle_+$
leading to the graph norm, i.e.
$\langle\phi_1,\phi_2\rangle_{+}
:=\langle(A^2+1)^{1/2}\phi_1,(A^2+1)^{1/2}\phi_2\rangle$.\par
Let $\N\subsetneq D(A)$ be a linear dense set which is 
closed with respect to the graph
norm. Being $\N$ closed we have $\HP=\N\oplus \N^\perp$ and we can then
consider the orthogonal projection $\pi:\HP\to \N^\perp$. 
More generally, since this gives advantages in concrete applications, 
we will consider a bounded linear map $\tau:\HP\to\fh$, where
$(\,\fh,[\cdot,\cdot]\,)$ is a Hilbert space. We suppose that $\tau$ is
surjective and that $\N:=$Ker$(\tau)$ is dense in $\H$. Note that, by
the surjectivity hypohtesis, $\fh\simeq \HP/\N
\simeq \N^\perp$, so that 
$\HP\simeq \N\oplus\fh$.\par
Denoting by $\rho(A)$ the resolvent set of $A$, for any $z\in\rho(A)$
we define the two bounded linear operators
$$
R(z):\H\to\HP\,,\qquad R(z):=(-A+z)^{-1}
$$
and 
$$
G(z):\H\to\fh\,,\qquad G(z):=\left(\tau R(\bar z)\right)^*\,.
$$
By \cite{[P3]}, Lemma 2.1 (see also \cite{[AKK]}, Theorem A.1, for an
analogous result), the denseness hypothesis on $\N$ is
equivalent to 
\begin{equation}
D(A)\cap\text{\rm Ran}(G(z))=\left\{0\right\}\,
\end{equation}
and, as an immediate consequence of the first resolvent 
identity for $R(z)$  (see \cite{[P1]}, Lemma 2.1), 
\begin{equation}
(z-w)\,R(w)G(z)=G(w)-G(z)\,.
\end{equation} 
These relations imply $$
D(A)\cap\text{\rm Ran}(G(w)+G(z))=\left\{0\right\}
$$
and 
$$
D(A)\supseteq\text{\rm Ran}(G(w)-G(z))\,.
$$
\begin{theorem} Defining 
$$R:= R(i)\,,\quad G:=G(-i)\,,\quad 
\G*:=\frac{1}{2}\,(G(i)+G(-i))\,,\quad
$$
one has 
$$
A^*_\N\phi=A\ph*+RG\zeta_\phi\,,
$$
$$
D(A^*_\N)=\left\{\,\phi\in\H\,:\,\phi=\ph*+\G* \zeta_\phi,\ \ph*\in
D(A),\ \zeta_\phi\in\fh\,\right\}\,.
$$
Defining
\begin{align*}
&\gamma_1:D(A^*_\N)\to\fh\,,\qquad\gamma_1\,\phi:=-\,\tau\ph*\,,\\
&\gamma_2:D(A^*_\N)\to\fh\,,\qquad\gamma_2\,\phi:=\zeta_\phi\,,
\end{align*}
the triple $\{\,\fh,\gamma_1,\gamma_2\,\}$ is a boundary triple 
for $A_\N^*$. The corresponding Weyl function of $A_\N$ is 
$$
\Gamma:\rho(A)\to\B(\fh)\,,\qquad \Gamma(z)=\tau(G_*-G(z))\,.
$$
\end{theorem}
\begin{proof} The form of $D(A^*_\N)$ was obtained in \cite{[P3]},
Theorems 3.4 and 4.1. By (3.2) 
\begin{equation}
RG=\frac{i}{2}\,(G(i)-G(-i))\,.
\end{equation}
Thus the action of $A^*_\N$ on its domain follows from \cite{[P3]}, 
Theorem 2.2.\par
Since $A$ is self-adjoint
$$
\langle A^*_\N\phi_*,\psi_*\rangle-
\langle \phi_*,A^*_\N \psi_*\rangle
=\langle A\phi_*,\psi_*\rangle-
\langle \phi_*,A\psi_*\rangle
=0\,.
$$
By (3.3)
\begin{equation}
G_*=-iRG+G\,.
\end{equation}
Therefore one has
\begin{align*}
&\langle A^*_\N G_*\zeta_\phi,G_*\zeta_\psi\rangle-
\langle G_*\zeta_\phi,A^*_\N G_*\zeta_\psi\rangle\\
=&\langle RG\zeta_\phi,G_*\zeta_\psi\rangle-
\langle G_*\zeta_\phi,R G\zeta_\psi\rangle\\
=&-2i\langle RG\zeta_\phi,RG\zeta_\psi\rangle
+\langle RG\zeta_\phi,G\zeta_\psi\rangle
-\langle G\zeta_\phi,R G\zeta_\psi\rangle\\
=&-2i\langle RG\zeta_\phi,RG\zeta_\psi\rangle
+\langle (-A-i)RG\zeta_\phi,RG\zeta_\psi\rangle\\
&-\langle (-A+i)RG\zeta_\phi,RG\zeta_\psi\rangle\\
=&\langle (2i+(-A-i)-(-A+i))RG\zeta_\phi,RG\zeta_\psi\rangle=0\,,
\end{align*}
\begin{align*}
&\langle A^*_\N\phi_*,G_*\zeta_\psi\rangle-
\langle \phi_*,A^*_\N G_*\zeta_\psi\rangle
=
\langle A\phi_*,G_*\zeta_\psi\rangle-
\langle \phi_*,RG\zeta_\psi\rangle\\
=&
-i\langle A\phi_*,RG\zeta_\psi\rangle+\langle A\phi_*,G\zeta_\psi\rangle-
\langle \phi_*,RG\zeta_\psi\rangle\\
=&
i\langle (-A-i)\phi_*,RG\zeta_\psi\rangle+\langle A\phi_*,G\zeta_\psi\rangle\\
=&
i\langle \phi_*,G\zeta_\psi\rangle+\langle A\phi_*,G\zeta_\psi\rangle
=
-\langle (-A+i)\phi_*,G\zeta_\psi\rangle\\
=&-[ G^*(-A+i)\phi_*,\zeta_\psi]
=-[\tau\phi_*,\zeta_\psi]\,
\end{align*}
and similarly 
$$
\langle A^*_\N G_*\zeta_\phi,\psi_*,\rangle-
\langle G_*\zeta_\phi,A^*_\N \psi_*\rangle
=[\zeta_\phi,\tau\psi_*]\,.
$$
In conclusion
$$
\langle A^*_\N \phi,\psi\rangle-
\langle \phi,A^*_\N \psi\rangle
=[-\tau\phi_*,\zeta_\psi]-[\zeta_\phi,-\tau\psi_*]\,.
$$
Since
$$
\K_z=\text{\rm Ker}(-A^*_\N+z)=\text{\rm Ran}(G(z))\,,
$$
for any 
$$\phi_z=G(z)\zeta\equiv(G(z)-G_*)\zeta+G_*\zeta$$ one has
$$
\Gamma(z)\gamma_2\phi_z=\tau(G_*-G(z))\zeta=-\tau(\phi_z)_*=\gamma_1\phi_z\,.
$$
\end{proof}
By the above theorem we can characterize all the singular
perturbations of $A$ by (abstract) boundary conditions of the kind 
$\Theta\,\zeta_\phi=\tau\phi_*$:
\begin{corollary} Any singular perturbations $\hat A$ of the self-adjoint
operator $A$ is of the kind 
$$
\hat A\phi=A\phi_*+RG\zeta_\phi\,,
$$ 
$$
D(\hat A)=\left\{\,\phi\in\H\,:\,\phi=\ph*+\G* \zeta_\phi,\ \ph*\in
D(A),\ \zeta_\phi\in\fh\,,\ \Theta\,\zeta_\phi=\tau\phi_*\right\}\,,
$$
where
$$
\tau\in\B(\H_+,\fh)\,,\quad \overline{\text{\rm Ker}(\tau)}=\H\,,\quad
\text{\rm Ran}(\tau)=\fh
$$
and $\Theta$ is a self-adjoint operator on the Hilbert space
$\fh$. Moreover 
$$
z\in\rho(\hat A)\cap\rho(A)\iff 0\in\rho(\Theta+\tau(G_*-G(z)))
$$ and

$$
(-\hat A+z)^{-1}=(-A+z)^{-1}+G(z)(\Theta+\tau(G_*-G(z)))^{-1}G(\bar
z)^*\,.
$$
\end{corollary}
\begin{remark} The above corollary was already obtained in \cite{[P3]} 
as a direct conseguence of the results contained in \cite{[P1]},
without making use of the theory of boundary triples.   
\end{remark}
We conclude by determining the eigenvectors (and their multiplicity)
of the singular perturbations of $A$. A similar result was obtained in
\cite{[BEKS]} in the case of form-bounded (hence weakly singular) 
perturbations. 
\begin{theorem} For any $\lambda\in\sigma_p(\hat A)\cap \rho(A)$, the map 
$\zeta\mapsto G(\lambda)\,\zeta$ is a bijection on
$\text{\rm Ker}(\Theta+\Gamma(\lambda))$ to 
$\text{\rm Ker}(-\hat A+\lambda)$.  
\end{theorem}
\begin{proof} 
At first let us note that 
$$
G(\lambda)\zeta\in
D(\hat A)\iff\tau(G(\lambda)-G_*)\zeta=\Theta\zeta
\iff(\Theta+\Gamma(\lambda))\zeta=0\,.
$$
We also note that, by (3.2) and (3.4),
$$
G(\lambda)-G_*=R(\lambda)(-\lambda G_*+RG)\,,
$$ 
i.e.
$$
(-A+\lambda)(G(\lambda)-G_*)=-\lambda G_*+RG\,.
$$
Suppose now that $\zeta\in\text{\rm Ker}(\Theta+\Gamma(\lambda))$. Then 
$G(\lambda)\zeta\in D(\hat A)$ and
\begin{align*}
&(-\hat A+\lambda)G(\lambda)\zeta=
(-\hat A+\lambda)((G(\lambda)-G_*)\zeta+G_*\zeta)\\
=&(-A+\lambda)(G(\lambda)-G_*)\zeta-RG\zeta+\lambda G_*\zeta\\
=&-\lambda G_*\zeta+RG\zeta-RG\zeta+\lambda G_*\zeta
=0\,.
\end{align*} 
Conversely suppose that $\phi\in \text{\rm ker}(-\hat A+\lambda)$. Then
\begin{align*}
0=&(-\hat A+\lambda)\phi=(-A+\lambda)\phi_*-RG\zeta_\phi+\lambda G_*\zeta_\phi\\
=&(-A+\lambda)\phi_*-(-A+\lambda)(G(\lambda)-G_*)\zeta_\phi\\
=&(-A+\lambda)(\phi_*+(G_*-G(\lambda))\zeta_\phi)\,.
\end{align*}
Since $\lambda\in\rho(A)$, this implies
$\phi_*+(G_*-G(\lambda))\zeta_\phi=0$, which is equivalent to 
$\phi=G(\lambda)\zeta_\phi$\,.
\end{proof}
\begin{remark} We know that the choice of a boundary triple is not
unique. Therefore a different (from the one given in Theorem
3.1) choice leads to a different parametrization of the family of the singular
perturbations of $A$. In particular the
resolvents will depend on a different Weyl function $\tilde\Gamma$. 
In \cite{[P1]} it was
obtained a Kre\u\i n-like formula which, in the terminology of the
present paper, gives singular perturbations of $A$ in terms of an
arbitrary choice of a Weyl function. Since, given any Weyl function 
$\tilde\Gamma$, one has that $\tilde\Gamma(z)-\Gamma(z)$ is a
$z$-independent self-adjoint operator (see \cite{[P1]}, remark 2.3),
if we denote by $\tilde A$ a singular perturbations of $A$
given by a Weyl function $\tilde\Gamma$, we have that 
$\tilde A=\hat A$, were 
$\hat A$ is the singular perturbation given in
Corollary 3.2 with $\Theta:=
\tilde\Theta+\tilde\Gamma(z)-\Gamma(z)$, $\tilde \Theta$
self-adjoint. Therefore Theorem 3.4 holds true 
independentently of the Weyl function 
(or, equivalently, of the boundary triple) one uses to 
describe the singular perturbations of $A$. 
\end{remark}



\end{document}